\documentclass[12pt]{article}
\usepackage{epsfig}
\def\const{{\mathrm{const}}}
\def\id{{\mathrm{id}}}
\def\R{\mathbf{R}}
\def\cross{\times}
\def\nu{\circ}
\def\C{{\mathbf{C}}}
\def\bC{{\mathbf{\overline C}}}
\def\mod{\mathrm{mod}\;}
\def\Ai{\mathrm{Ai}}
\title{Nevanlinna functions with real zeros}
\author{A. Eremenko\thanks{Supported by NSF grants
DMS-0100512 and DMS-0244547.}\, 
and S. Merenkov\thanks{
Supported by NSF grants DMS-0400636 and DMS-0244547.}
}
\date{January 2005}
\begin{document}
\maketitle

\section{Introduction}

Consider a differential equation
\begin{equation}
\label{DE}
w''+Pw=0,
\end{equation}
where $P$ is a polynomial of the independent variable. Every solution
$w$ of this equation is an entire function.
We are interested in solutions $w$ whose all roots are real.
If (\ref{DE}) has two linearly independent solutions with
this property then $\deg P=0$, see \cite{G2,HR}.
Here we study equations (\ref{DE}) that have
at least one solution with all roots real.

The question of describing equations (\ref{DE}) with this property
was proposed by S. Hellerstein
and J. Rossi in \cite[Probl. 2.71]{H}. According to \cite{G1},
up to trivial changes of variables,
only the following four examples were known until recently. 
\begin{itemize}
\item
$\deg P=0.$ If $k$ is real, all solutions of $w''+k^2w=0$
are trigonometric functions.
\item
$\deg P=1$. The Airy equation $w''-zw=0$,
has a solution $\Ai(z)$ whose
roots lie on the negative ray.
\item
$\deg P$ is even and $w=p\exp q$,
where $p$ and $q$ are polynomials, and all
roots of $p$ are real. In this case, the set
of roots of $w$ is evidently
finite.
For example, if $P(z)=1-z^2+2n,$ where $n$ is a positive integer,
then
the equation (\ref{DE}) has solutions $w=H_n(z)\exp(-z^2/2),$
where $H_n$ are Hermite's polynomials whose roots are all real.

\item
$P(z)=az^4+bz^2-c$. Gundersen \cite{G1} proved that 
for every $a>0$ and $b\geq 0$ one can find
an infinite set of real numbers $c$, such that some solution
of $(\ref{DE})$ 
has infinitely many roots,
almost all of them real. When $b=0$ this result was earlier obtained
by Titchmarsh \cite[p. 172]{T}.
\end{itemize}

Here and in what follows ``almost all'' means ``all except finitely many''.
Recently Kwang C. Shin \cite{shin} proved a similar result for
a degree $3$ polynomial:
\begin{itemize}

\item
For every real $a$ and $b\leq 0$ there exist an infinite set
of positive numbers $c$ such that the equation
$w''+(z^3+az^2+bz-c)w=0$ has a solution with infinitely many roots,
almost all of them real.
\end{itemize}

On the other hand, Gundersen \cite{G2} proved the following theorems:
\vspace{.1in}

\noindent
{\bf Theorem A} {\em If $d=\deg P\equiv 2\; (\mod\, 4),$
and $w$ is a solution of $(\ref{DE})$, with almost all roots real,
then $w$
has only
finitely many roots.} 
\vspace{.1in}

\noindent
{\bf Theorem B} {\em If $(\ref{DE})$ possesses a solution $w$
with infinitely many real zeros,
then $P$ is a real polynomial, and $w$
is proportional to a real function.}
\vspace{.1in}

We also mention a result of Rossi and Wang \cite{RW} that
if (\ref{DE}) has a solution with infinitely many roots, all of
them real, then at least half of all roots of $P$ are non-real.
In view of Theorem B we restrict from now on to the case
of real polynomials $P$ in (\ref{DE}). Our results are: 
\vspace{.1in}

\noindent
{\bf Theorem 1} {\em For every $d$,
there exist $w$ satisfying $(\ref{DE})$ with $\deg P=d$ and such that
all roots of $w$ are real.
For every positive $d$ divisible by $4$, there exist $w$ with infinitely
many roots, all of them real, as well as $w$ with any prescribed finite
number of roots, all of them real.}
\vspace{.1in}

\noindent
{\bf Theorem 2} {\em Let $w$ be a solution of
the equation $(\ref{DE})$ 
whose all roots are real, and $d=\deg P$. Then:
\newline
(a) For $d\equiv 0\,(\mod 4)$ the set of roots of
$w$ is either finite or unbounded from above and below (as a subset of the real axis). 
\newline
(b) For odd $d$ the roots of $w$ lie on a ray, and there
are infinitely many of them.}
\vspace{.1in}

Theorem 2 can be generalized to the case that almost all roots
of $w$ are real.

In comparison with the existence results of Titchmarsh, Gundersen and Shin
mentioned above, our Theorem 1 gives more precise information on the
zeros $w$: they are all real. On the other hand we
can tell less about the polynomial $P$.

Our proofs are based on a geometric
characterization of meromorphic functions of
the form
$f=w_1/w_2$, where $w_1$ and $w_2$ are linearly independent solutions of
(\ref{DE}), due to F. and R. Nevanlinna \cite{Fne,RNe1}, which will be
explained in the next section.

\section{A class of meromorphic functions}

We associate with (\ref{DE}) another differential equation
\begin{equation}
\label{SDE}
\frac{f'''}{f'}-\frac{3}{2}\left(\frac{f''}{f'}\right)^2=2P.
\end{equation}
The expression in the left hand side of (\ref{SDE}) is called the
{\em Schwarzian derivative} of $f$. The following well-known fact
is proved by simple formal computation.
\vspace{.1in}

\noindent
{\bf Proposition 1} {\em The relation $f=w_1/w_2$ gives a
bijective correspondence
between solutions $f$ of $(\ref{SDE})$ and classes of proportionality of
pairs $(w_1,w_2)$ of
linearly independent solutions of $(\ref{DE})$.}
\vspace{.1in}

Thus when $P$ is a polynomial, all solutions of (\ref{SDE}) are
meromorphic in the complex plane, and they are all obtained from each other
by post-composition with a fractional-linear transformation.
We call solutions of equations (\ref{SDE}) with polynomial right hand side
{\em Nevanlinna functions}.
Equation (\ref{SDE}) has a real solution if and only if $P$ is real.

It is easy to see that meromorphic functions $f$ satisfying (\ref{SDE}) are
local homeomorphisms. In other words, $f'(z)\neq 0$ and all poles are
simple.

F. and R. Nevanlinna gave a topological characterization
of all meromorphic functions $f$
which may occur as solutions of (\ref{SDE}). To formulate their result
we recall several definitions.

A {\em surface} is a connected Hausdorff topological manifold of
dimension $2$ with countable base.

A continuous map $\pi:X\to Y$ of surfaces is called
{\em topologically holomorphic} if it is open and discrete.
According to a theorem of Stoilov \cite{S} this is equivalent to
the following property. For every $x\in X$, there is a positive integer $k$ and
complex local coordinates $z$ and $w$ in neighborhoods of $x$ and $\pi(x)$,
such that $z(x)=0$ and the map $\pi$
has the form $w=z^k$ in these coordinates. The integer $k=k(x)$ is called
the
{\em local degree} of $\pi$ at the point $x$. So $\pi$ is
a local homeomorphism if and only if $n=1$ for every $x\in X$.

A pair $(X,\pi)$ where $X$ is a surface and $\pi:X\to\bC$ a topologically
holomorphic map 
is called {\em a surface spread over the sphere} 
(\"Uberlagerungsfl\"ache in German).
Two such pairs $(X_1,\pi_1)$
and $(X_2,\pi_2)$ are considered
equivalent if there is a homeomorphism $h:X_1\to X_2$ such that  
$\pi_2=\pi_1\circ h$.
So, strictly speaking, a surface spread over the sphere
is an equivalence class of such pairs.

If $f:D(R)\to\bC$ is a meromorphic function in some
disc $D(R)=\{ z:|z|<R\},\; R\leq\infty$ then $(D(R),f)$ defines a
surface spread over the sphere.
We will call the equivalence class of $(D(R),f)$ the
{\em surface
associated with $f$}. It is the same as the Riemann surface of $f^{-1}$,
as it is defined on \cite[p. 288]{Ah}, completed with algebraic branch
points as in \cite[p. 300]{Ah}.

In the opposite direction, suppose that $(X,\pi)$ is a surface spread over
the sphere. Then there exists a unique conformal structure on $X$ which makes
$\pi$ holomorphic. If $X$ is open and simply connected, the Uniformization
Theorem says that there exists a conformal homeomorphism
$\phi:D(R)\to X$, where $R=1$ or $\infty$.
This $\phi$ is defined up to a conformal automorphism of $D(R)$.
The function $f=\pi\circ\phi$ is meromorphic in $D(R)$, and $(X,\pi)$
is (a representative of) the surface associated with $f$.

If $R=\infty$ we say that $(X,\pi)$ is of {\em parabolic type}.

We consider surfaces spread over the sphere $(X,\pi)$ where
$X$ is open and simply connected\footnote{That is homeomorphic to the plane.},
$\pi$ a local homeomorphism, and subject to additional conditions below.

Suppose that for some finite set $A\subset\bC$ the
restriction
\begin{equation}
\label{S}
\pi:X\backslash\pi^{-1}(A)\to\bC\backslash A\quad\mbox{is a covering map}.
\end{equation}
Fix an  open topological disc 
$D\subset \bC$ containing exactly one point $a$ of the set $A$.
If $V$ is a connected component
of $\pi^{-1}(D\backslash \{ a\})$ then
the restriction
\begin{equation}
\label{V}
\pi:V\to D\backslash \{ a\}
\end{equation}
is a covering of a ring, and its degree $k$ 
does not depend on the choice of the disc $D$.
The following cases are possible.
\newline
a) $k=\infty$. Then (\ref{V}) is a universal covering, $V$ is simply
connected and its boundary consists of a single simple curve in $X$
tending to ``infinity'' in both directions. In this case we say that $V$
defines a {\em logarithmic singularity} over $a$.
Notice that the number of logarithmic singularities over $a$
is independent of the choice of $D$.
\newline
b) If $k<\infty$ and there exists a point $x\in X$ such that
${\widetilde{V}}=V\cup\{ x\}$ is an open topological disc, then
$\pi:{\widetilde{V}}\to D$ is a ramified covering and
has local degree $k$ at $x$. As we assume that $\pi$ is a local
homeomorphism, only $k=1$ is possible, so 
$\pi:{\widetilde{V}}\to D$ is a homeomorphism.
\newline
c) If $k<\infty$ but there is no point $x\in X$ such that $V\cup\{ x\}$
is an open disc, then we can add to $X$ such an ``ideal point'' and define
the topology on ${\widetilde{X}}=X\cup\{ x\}$ so that it remains a surface.
Evidently ${\widetilde{X}}$ is a sphere, and our local homeomorphism
extends to a topologically holomorphic map between spheres whose
local degree equals one everywhere except possibly one point.
It easily follows that the local degree equals one everywhere,
the extended map is a homeomorphism. This implies that our original
map $\pi$ is an embedding. 

So in any case the degree of the map (\ref{V})
can be only $1$ or $\infty$. 
\vspace{.1in}

\noindent
{\bf Definition}
We say that $(X,\pi)$ is an {\em N-surface} if $X$ is open and simply
connected,
$\pi$ is a local homeomorphism, condition (\ref{S}) is
satisfied, and there
are only finitely many logarithmic singularities. 
\vspace{.1in}

A simple topological argument shows that aside from the case when
$\pi$ is an embedding, the number of logarithmic singularities is at least two.
All cases with two logarithmic singularities can be reduced by
a fractional-linear transformation of $\bC$ to the case
$\exp:\C\to\bC$.

The name N-surface is chosen in honor of F. and R. Nevanlinna.
The complete official name of this object would be ``An open simply
connected surface spread over the sphere without algebraic branch points
and with finitely many logarithmic singularities''.
\vspace{.1in}

\noindent
{\bf Theorem C} (i) {\em Every N-surface is of parabolic type, that is its
associated functions are meromorphic in the plane $\C$.}
\newline
(ii) {\em  If an N-surface $(X,\pi)$ has $n\geq 2$ logarithmic singularities then the
associated meromorphic functions $f=\pi\circ\phi$ satisfy a differential
equation $(\ref{SDE})$ in which $\deg P=n-2.$}
\newline 
(iii) {\em For every polynomial $P$, every solution $f$ of $(\ref{SDE})$ is
a meromorphic function in the plane whose associated surface
is an N-surface with $n=\deg P+2$ logarithmic singularities.}
\vspace{.1in}

Two meromorphic functions $f_1$ and $f_2$ are called
equivalent if $f_1(z)=f_2(az+b)$ with $a\neq 0$.
Theorem C establishes a
bijective correspondence between equivalence 
classes of Nevanlinna functions and N-surfaces.

In this paper we use only statements (i) and (ii) of Theorem C.
The connection between N-surfaces and differential equations
was apparently discovered by F. Nevanlinna who proved (iii) in 
\cite{Fne}.
Statements (i) and (ii) were proved for the first time by
R. Nevanlinna in \cite{RNe1}.
Then Ahlfors
\cite{Ah1} gave an alternative proof based on completely different ideas.
A modern version of this second proof uses quasiconformal mappings
\cite{GO}. This modern proof is reproduced in \cite{DV}.
All these authors were primarily interested in the theory of
meromorphic functions, and used differential equations as a tool.
Apparently, the only application of Theorem C to differential equations
is due to Sibuya \cite{Sib} who deduced from it the existence of
equations (\ref{DE}) with prescribed Stokes multipliers.

In view of Theorem C, to obtain our results, we only need to
single out those N-surfaces that are associated with
real meromorphic functions with real zeros.

\section{Speiser graphs}

We recall a classical tool for explicit construction and visualization
of N-surfaces. 
It actually applies to all surfaces spread over the
sphere that satisfy (\ref{S}).
First we 
suppose that a surface spread over the sphere $(X,\pi)$ satisfying (\ref{S}) is given, and $A=\{ a_1,\ldots,a_q\}$
is the set in (\ref{S}). We call elements of $A$ {\em base points}.
Consider a {\em  base curve} that is an oriented Jordan 
curve $\Gamma$ passing
through $a_1,\ldots,a_q$. Choosing a base curve defines a cyclic order
on $A$, and we assume that the enumeration is consistent with this
cyclic order and interpret the subscripts as remainders modulo $q$.

The base curve $\Gamma$ divides the Riemann sphere $\bC$ into two
regions which we denote $D^\cross$ and $D^\circ$, so that
when $\Gamma$ is traced according to its orientation, the region $D^\cross$ is
on the left. We choose points
$\times\in D^\cross$ and $\circ\in D^\circ$, and connect these
two points by $q$ disjoint simple arcs $L_j$ so that each $L_j$ intersects
$\Gamma$ at exactly one point, and this point belongs to the arc
$(a_j,a_{j+1})\subset\Gamma$.
We obtain an {\em embedded graph} $L\subset\bC$ having two
vertices $\cross$ and $\circ$ and $q$ edges $L_j$.
This embedded graph defines a cell decomposition of the sphere,
whose $2$-cells (faces) are components of the complement of $L$,
$1$-cells (edges) are the open arcs $L_j$ and $0$-cells (vertices)
are the points $\cross$ and $\circ$.
Each face contains exactly one base point. 

The preimage of this cell decomposition under $\pi$ is a cell
decomposition of $X$, because as we saw in the previous section,
each component of the preimage of a cell is a cell of the same dimension.
The $1$-skeleton $S=\pi^{-1}(L)\subset X$
is a connected
properly embedded graph in $X$.
As $S$ completely defines the cell decomposition, we will permit ourselves to
follow the tradition and speak of this graph instead of the cell decomposition,
and use such expressions as
``faces of $S$'' meaning the faces of the cell decomposition.

We label vertices of $S$ by $\cross$ and $\circ$, according to their images
under $\pi$, and similarly label the faces by the base points $a_j$.
Our labeled graph $S$ (or more precisely, the labeled cell decomposition)
has the following properties:
\vspace{.1in}

\noindent
1. Every edge connects a $\cross$-vertex
to a $\circ$ vertex.
\vspace{.1in}

\noindent
2. Every vertex belongs to the boundaries
of exactly $q$ faces having all $q$ different labels.
\vspace{.1in}

\noindent
3. The face labels have cyclic order
$a_1,\ldots,a_q$ anticlockwise around each $\cross$-vertex, and the opposite
cyclic order around each $\circ$-vertex.
\vspace{.1in}

The labeled graph $S$ is called the {\em Speiser
graph} or the {\em line complex} of the surface spread over the sphere
$(X,\pi)$.

A face of $S$ is called {\em bounded} if its boundary consists of finitely
many edges and vertices. It follows from property 1 that the numbers
of edges and vertices on the boundary of a bounded
face are equal and even.
If $a$ is a base point, all solutions of the equation $\pi(x)=a$
belong to the bounded faces labeled by $a$,
and each such face contains exactly
one solution. 
If $k$ is the local degree of $\pi$ at this point $x$ then the face
is a $2k$-gon, that is bounded by $2k$ edges and $2k$ vertices. 

Suppose now that $X$ is a surface and a labeled cell decomposition of $X$
with $1$-skeleton $S$
is given such that 
1, 2 and 3 hold.
If we choose a set $A\subset\bC$ of $q$ points and a
curve $\Gamma\subset\bC$ passing through the
points of the set $A$,
and define 
$L$ as above, then there exists a
topologically holomorphic map $\pi$ such that 
$S=\pi^{-1}(L)$. This map $\pi$
is unique up to pre-composition with a homeomorphism
of $X$.
A verification of this statement is contained in
\cite{GO}. We recall the construction.

The labels of faces define labels of edges: an edge is labeled by $j$
if it belongs to the common boundary of faces with labels $a_j$ and $a_{j+1}$.
This defines a map of the $1$-skeleton to the $1$-skeleton $L$
of the cell decomposition of the sphere: each edge of $S$ labeled by $j$
is mapped onto $L_j$ homeomorphically, and such that the orientation
is consistent with the vertex
labeling. It is easy to see that this map is a covering $S\to L$.
The boundary of each face covers a topological circle formed by two
adjacent
edges and two vertices of $L$. Such map extends to a ramified covering
between faces, ramified only over the base points (unramified for
$2$-gonal or unbounded faces).

For a given cyclically ordered set $A$ and a surface $X$,
the correspondence between topologically holomorphic maps $\pi$
and Speiser graphs
is not canonical: it depends on the choice of the base curve $\Gamma$.
(It is the isotopy class of $\Gamma$ with fixed set $A$ that matters).

It is easy to single out those Speiser graphs that correspond to N-surfaces:
the ambient surface $X$ is open and simply connected, and two
additional properties hold:
\vspace{.1in}

\noindent
4. Each face has either two or infinitely many boundary edges.
\vspace{.1in}

\noindent
5. The set of unbounded faces is finite. 
\vspace{.1in}

Property 4 corresponds to the assumption that $\pi$ is a local homeomorphism,
and property 5 follows from the fact that unbounded faces correspond to
``logarithmic singularities'' that is to the components $V$ in (\ref{V})
where the covering has infinite degree.

So we have
\vspace{.1in}

\noindent
{\bf Proposition 2} {\em Let $(X,\pi)$ be an N-surface, and $a$ a basis point.
Then each solution of the equation $\pi(x)=a$ is contained in a face
which is a $2$-gon, and is labeled by $a$. Each such face contains exactly
one solution of this equation.}
\vspace{.1in}

For each Speiser graph $S$ corresponding to an N-surface, we 
define a new graph $T(S)$ with the same vertices:
two vertices are connected by a single edge in $T$ if they
are connected by at least one edge in $S$.
Thus $T(S)$ is obtained from $S$ by dropping multiple edges.
Property 4 of $S$ implies that $T$ is a tree. Faces of $T(S)$ are
exactly the unbounded faces of $S$.
We assume that vertices and unbounded faces inherit their labels from
$S$. 

The tree $T(S)$ is properly embedded in $X$. Each vertex 
has at least two and at most $q$ adjacent edges in $T(S)$
and the cyclic order of face labels around the $\cross$-vertices
of $T(S)$ is the same as in $S$.

Suppose that $S$ has more than two unbounded faces. Then
the tree $T(S)$ has $n$ maximal (by inclusion) infinite ``branches'' having
all vertices of degree $2$, of the
form 
$$-\circ-\cross-\circ-\cross-\circ-\ldots\quad\mbox{or}\quad
-\cross-\circ-\cross-\circ-\cross-\ldots$$
where $n$ is the number of logarithmic singularities of $(X,\pi)$.
Such branches are called {\em logarithmic ends}.
A tree $T(S)$ is the union of its logarithmic ends and a finite
subtree.

\section{Symmetric Speiser graphs}

A {\em symmetric} surface spread over the sphere is defined as
a triple $(X,\pi,s)$, where 
$s:X\to X$ is a homeomorphism such that $s\circ s=\id$
and $\pi\circ s(x)=\overline{\pi(x)}$,
and the bar denotes complex conjugation. It is clear that such $s$ is
an anticonformal homeomorphism of $X$. If $(X,\pi)$ is
of parabolic type, and $\phi:\C\to X$ a conformal homeomorphism
then $\phi^{-1}\circ s\circ\phi$ is an anticonformal involution of the
complex plane. Each such involution is conjugate to $z\mapsto\overline{z}$
by a conformal automorphism of $\C$. So for a symmetric surface
spread over the sphere there exists a uniformizing map $\phi$ with
the property $\phi(\overline{z})=s\circ\phi(z),\, z\in\C$.
The set of fixed points of $s$ is called the {\em axis} (of symmetry);
it is the image of the real line under $\phi$.

If $f$ is a real function meromorphic in $\C$ then its associated
surface has a natural involution which makes it a symmetric
surface spread over the sphere. In the opposite direction,
to a symmetric surface spread over the sphere of parabolic type,
a real meromorphic function is associated. 

It is clear that the set $A=\{ a_1,\ldots,a_q\}$ of basis points 
of a symmetric N-surface is invariant under complex conjugation.
Suppose for a moment that at most two of the points $a_1,\ldots,a_q$ are real.
Then there exists a base curve $\Gamma$ passing through $a_1,\ldots,a_q$
which is symmetric with respect to complex conjugation.
Choosing the $\cross$ and $\circ$ points on the real axis we
can perform the construction of the Speiser graph symmetrically.
The resulting graph $S$ and the tree $T(S)$ will be {\em symmetric} in
the following sense. The involution $s$ will send each vertex to
a vertex with the same label, each edge to an edge and each face to a face
with complex conjugate label.

In the general case, that more than two basis points are allowed
on the real line, one has to modify a little the definition of the Speiser 
graph. Let $a$ be a real base point of a symmetric N-surface $(X,\pi,s)$.
Consider an open (round) disc $D\subset\bC$ centered at $a$ and not
containing other base points. Let $V$ be a component of $\pi^{-1}(D)$
which defines a logarithmic singularity (see Section 2).
Then one of the following: either $V$ is invariant under $s$, or $s(V)=V'$
where $V'$ is another component of $\pi^{-1}(D)$, disjoint from $V$. In
the latter case, $V$ is disjoint from the symmetry axis.
In the former case, we will call the logarithmic singularity {\em real}.

We claim that there are at most two real logarithmic singularities.
Indeed, for a real logarithmic singularity, the intersection of
$V$ with the symmetry axis consists of a ``ray'', and there cannot be more
than two disjoint ``rays'' on the symmetry axis.

The symmetry axis divides $X$ into two ``halfplanes'', and each non-real
logarithmic singularity (more precisely, its defining region $V$)
belongs to one of these ``halfplanes''. Thus
the non-real logarithmic singularities are split into two classes, say
$C_+$ and $C_-$
according to the ``halfplane'' they belong, and the regions $V$ and $V'$
always belong to different classes.

The real logarithmic singularities lie over at most two
basis points.

Let $a$ be a real basis point
such that there are no
real logarithmic singularity over $a$. Consider a homeomorphism $\eta_+$ of
the Riemann sphere which is identical outside $D$ and sends the point $a$
to the point $a+i\epsilon$, where $\epsilon>0$ is so small that $a+i\epsilon\in
D$. Let $\eta_-(z)=\overline{\eta_+(\overline{z})}$.
We deform our map $\pi$ in the following way:
$$\pi^*(x)=\left\{\begin{array}{ll}\eta_+\circ\pi(x),\; x\in V, &V\in C_+,\\
                                   \eta_-\circ\pi(x),\; x\in V, &V\in C_-,\\
                                   \pi(x)  &\mbox{otherwise.}
\end{array}\right.
$$
Evidently, the new N-surface is symmetric.
Let $E$ be the set projections of real logarithmic singularities,
${\mathrm{card}\; E}\leq 2$.
Performing the deformation described above for all real base points
except two of them, $a'$ and $a''$, such that $E\subset\{ a',a''\}$ we
obtain a new symmetric $N$-surface which has the property
that only two basis points are real. So a symmetric base curve
can be chosen and a symmetric Speiser graph constructed. 
The Speiser graph of this deformed surface
does not depend on $\epsilon$
as soon as $\epsilon$ is positive and small enough, and we call
it a {\em symmetric Speiser graph of} $(X,\pi)$. The number of
basis points of  $(X,\pi^*)$ is larger than that of
$(X,\pi)$; to preserve properties 2 and 3, we can use two
different
labels, say $a^+$ and $a^-$ for a basis point $a$ as above.
The faces of $S$ over $a+i\epsilon$ are labeled by $a^+$, those
over $a-i\epsilon$ 
by $a^-$.

Symmetric Speiser graphs have all the properties 1-5 listed above, and
in addition, they are preserved by an orientation-reversing involution
of the ambient surface.
This includes the action of $s$ on labels
if we consider
$a^+$ and $a^-$ as complex conjugate.

Given a symmetric Speiser graph 
one can construct
a symmetric surface spread over the sphere corresponding to this Speiser
graph. First we replace all labels $a^+$ by $a+i\epsilon$ and $a^-$ by 
$a-i\epsilon$, then choose a symmetric basis curve passing through the
new basis points, points $\cross$ and $\circ$ on the real axis
and a symmetric graph $L$, and perform all
construction preserving symmetry.
Then we apply the inverse of the deformation described above to place the basis points
in their original position. 

We will need two simple properties of symmetric Speiser graphs $S$ and 
trees $T(S)$.
\vspace{.1in}

\noindent
A. If a logarithmic end intersects the axis then it is contained in
the axis.
\vspace{.1in}

\noindent
B. Every edge either belongs to the axis or is disjoint from it.
\vspace{.1in}

Consider a symmetric N-surface satisfying the following
\vspace{.1in}

\noindent
{\bf Assumptions}
{\em The number of logarithmic ends is at least $3$,
zero is a basis point, there is at most one
real logarithmic singularity not lying over zero, and the symmetric
Speiser graph does not have labels $0^+$ or $0^-$. }
\vspace{.1in}

Comments. The assumption that there are at least three logarithmic ends
excludes only the trivial cases when $P=\const$.
The assumption that $0$ is a basis point does not restrict generality
because an extra basis point can be always added.
The third assumption excludes the cases when the number of real zeros
is finite (if there are two non-zero real logarithmic singularities,
then $f$ has non-zero limits along the real axis as $x\to+\infty$
and $x\to-\infty$,
so the set of real zeros is finite). If this third assumption is
satisfied, we can always construct the symmetric
Speiser graph in such a way that it does not contain
labels $0^+$ and $0^-$.
\vspace{.1in}

\noindent
{\bf Proposition 3}. 
{\em Let $S\subset\C$ be a symmetric Speiser graph corresponding to 
a symmetric N-surface, $f$ an associated real Nevanlinna function,
and the above Assumptions are satisfied.
Then all zeros of $f$ are real if and only if $S$ has the following property:
each vertex belongs either to the axis of symmetry, or to the boundary
of an unbounded face labeled by $0$.}
\vspace{.1in}

{\em Proof}. Suppose that all zeros of $f$ are real.
This means that all $\pi$-preimages of $0$ lie on the axis of symmetry.
By Proposition 2, these preimages are in bijective
correspondence with $2$-gonal faces $F$ labeled by zero.
We have for such faces $F\cap s(F)\neq\emptyset$,
and thus by the symmetry of the graph,
$F=s(F)$. It follows that both vertices on the boundary of
$F$ belong to the axis.

Every vertex belongs to the boundary of some face
labeled by $0$. If this face is bounded we conclude from the above
that the vertex lies
on the symmetry axis. This proves necessity of the condition of Proposition 3.

Now we suppose that each vertex belongs either to the symmetry axis or to
the boundary of an unbounded face labeled by zero. As every vertex
belongs to the boundary of {\em only one} face labeled by $0$, we conclude
that every $2$-gonal face labeled by zero has one boundary vertex
on the axis, and thus 
its other boundary vertex also belongs to the axis.
We conclude that this face is symmetric, and thus the $\pi$-preimage of $0$ contained in this face belongs to the axis.
\vspace{.1in}

{\em Proof of Theorem 2}. The cases $P=0$ and $d=0$ are trivial,
so we assume
that $d\geq 1$. If there are two real logarithmic singularities
over non-zero points then $d$ is even and $f$ has finitely
many zeros, so there is nothing to prove.
Thus we suppose from now on that the Assumptions stated above
are satisfied.

If the number $n=\deg P+2$ of logarithmic ends is even,
then either none or two of them belong to the axis of symmetry.

If none of the logarithmic ends belongs to the axis then
$f$ has finitely many zeros.

Now suppose that there are two logarithmic ends on the axis
of symmetry.
Consider one of them.
Let $a$ and $b$ be the labels of the two unbounded
faces adjacent to it.
Then $a\neq b$ and $a=\overline{b}$, so neither 
$a$ nor $b$ can be $0$.
This implies that all vertices on this logarithmic end belong to the
boundaries of $2$-gonal faces whose labels are $0$, and
thus we obtain an infinite sequence of real zeros.
As there are two logarithmic ends on the axis of symmetry,
the sequence of zeros is unbounded from above and below.
This proves (a).

That for $d=4k, k\geq 1$ both cases actually occur
is demonstrated by Figures 1 and 2.

If $n$ is odd then exactly one logarithmic
end belongs to the axis of symmetry. 
The other end of the symmetry axis is contained in an infinite face
and bisects it (by symmetry).
It follows that $\pi$ has a limit on this other
end of the symmetry axis.
If the limit is non-zero, (b) immediately follows.
If the limit is zero, we notice that this end of the axis belongs to
a neighborhood $V$ of a real logarithmic singularity over $0$,
and again statement (b) follows.
\vspace{.1in}

{\em Proof of Theorem 1}. It is enough to display a Speiser graph for each case.
For simplicity we only show in figures 1-5 the trees $T(S)$ on the left of each picture
and the basis curves with basis points
on the right. It can be easily seen that each of our trees has unique
extension to a symmetric Speiser graph $S$ that has the property
described in Proposition 3.
\vspace{.1in}

{\em Remarks}. Suppose that all basis points of a symmetric N-surface
are real. Then we can construct
another kind of Speiser graph which we call {\em almost symmetric},
without using the perturbation procedure described above. Namely, take the
real axis as the basis curve, and choose $\cross$ and $\circ$ at the
points $\pm i$. The corresponding Speiser graph is preserved by the
involution,
except that the vertex labels are now interchanged. It is easy to see
that in such almost symmetric graph there can be no vertices on the axis,
and {\em exactly one edge} of $T(S)$ intersects the axis. We conclude
that

{\em A real Nevanlinna function with at least three logarithmic
singularities and only real asymptotic values 
can have at most one
real zero}.
\newpage

\begin{center}
\begin{picture}(0,0)%
\includegraphics{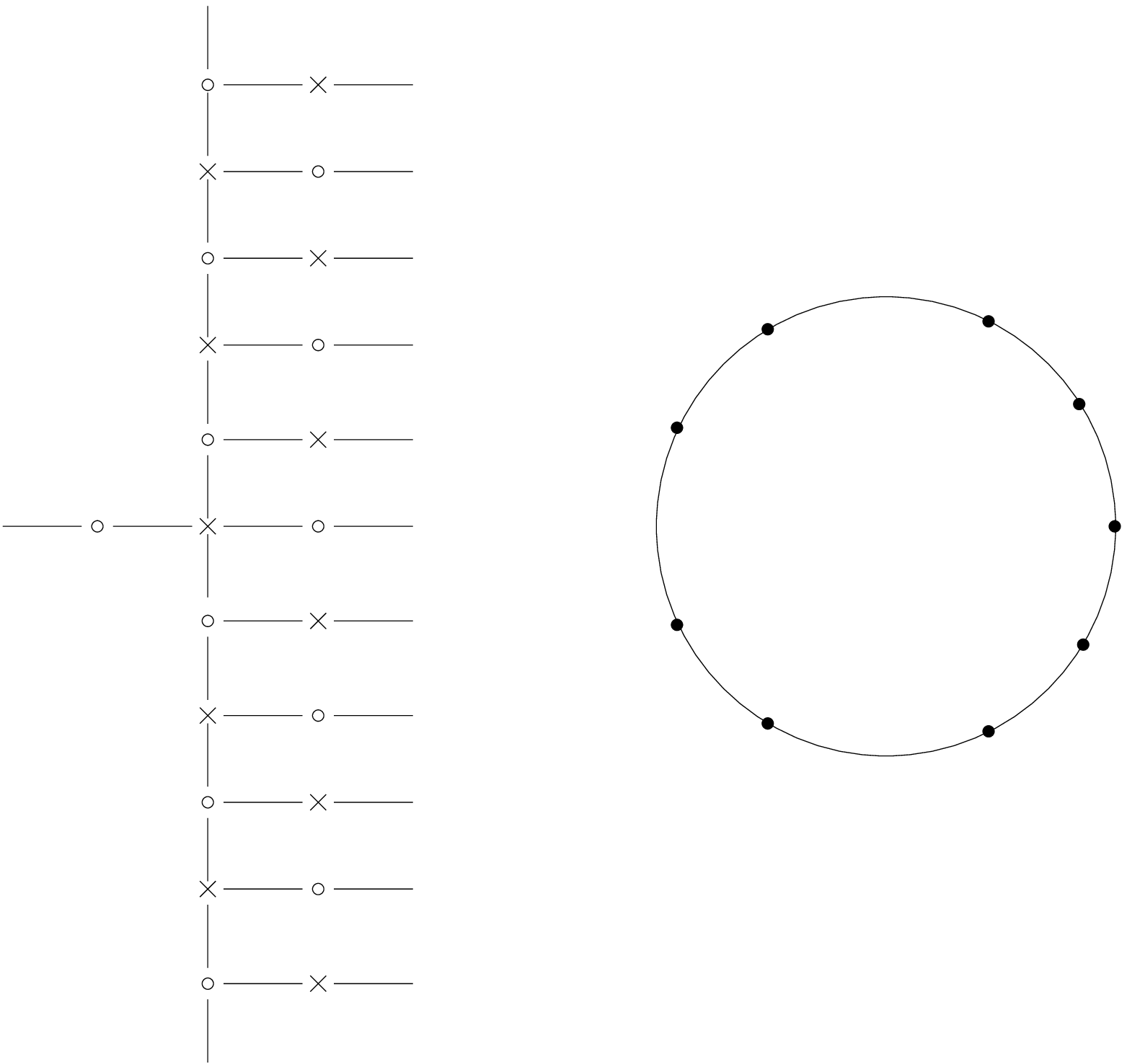}%
\end{picture}%
\setlength{\unitlength}{2763sp}%
\begingroup\makeatletter\ifx\SetFigFont\undefined%
\gdef\SetFigFont#1#2#3#4#5{%
  \reset@font\fontsize{#1}{#2pt}%
  \fontfamily{#3}\fontseries{#4}\fontshape{#5}%
  \selectfont}%
\fi\endgroup%
\begin{picture}(10737,10074)(1039,-9373)
\put(3301,-3961){\makebox(0,0)[lb]{\smash{\SetFigFont{12}{14.4}{\rmdefault}{\mddefault}{\updefault}
$c$
}}}
\put(3376,-4861){\makebox(0,0)[lb]{\smash{\SetFigFont{12}{14.4}{\rmdefault}{\mddefault}{\updefault}
$\overline{c}$
}}}
\put(3301,-3061){\makebox(0,0)[lb]{\smash{\SetFigFont{12}{14.4}{\rmdefault}{\mddefault}{\updefault}
$0$
}}}
\put(3301,-1486){\makebox(0,0)[lb]{\smash{\SetFigFont{12}{14.4}{\rmdefault}{\mddefault}{\updefault}
$0$
}}}
\put(3226,164){\makebox(0,0)[lb]{\smash{\SetFigFont{12}{14.4}{\rmdefault}{\mddefault}{\updefault}
$0$
}}}
\put(3376,-5836){\makebox(0,0)[lb]{\smash{\SetFigFont{12}{14.4}{\rmdefault}{\mddefault}{\updefault}
$0$
}}}
\put(3301,-7411){\makebox(0,0)[lb]{\smash{\SetFigFont{12}{14.4}{\rmdefault}{\mddefault}{\updefault}
$0$
}}}
\put(3301,-9136){\makebox(0,0)[lb]{\smash{\SetFigFont{12}{14.4}{\rmdefault}{\mddefault}{\updefault}
$0$
}}}
\put(3301,-2311){\makebox(0,0)[lb]{\smash{\SetFigFont{12}{14.4}{\rmdefault}{\mddefault}{\updefault}
$b_1$
}}}
\put(3301,-661){\makebox(0,0)[lb]{\smash{\SetFigFont{12}{14.4}{\rmdefault}{\mddefault}{\updefault}
$b_2$
}}}
\put(6976,-3361){\makebox(0,0)[lb]{\smash{\SetFigFont{12}{14.4}{\rmdefault}{\mddefault}{\updefault}
$a$
}}}
\put(11401,-5611){\makebox(0,0)[lb]{\smash{\SetFigFont{12}{14.4}{\rmdefault}{\mddefault}{\updefault}
$\overline{c}$
}}}
\put(10426,-6511){\makebox(0,0)[lb]{\smash{\SetFigFont{12}{14.4}{\rmdefault}{\mddefault}{\updefault}
$\overline{b_1}$
}}}
\put(3301,-8236){\makebox(0,0)[lb]{\smash{\SetFigFont{12}{14.4}{\rmdefault}{\mddefault}{\updefault}
$\overline{b_2}$
}}}
\put(3301,-6586){\makebox(0,0)[lb]{\smash{\SetFigFont{12}{14.4}{\rmdefault}{\mddefault}{\updefault}
$\overline{b_1}$
}}}
\put(1276,-5461){\makebox(0,0)[lb]{\smash{\SetFigFont{12}{14.4}{\rmdefault}{\mddefault}{\updefault}
$\overline{a}$
}}}
\put(1351,-3061){\makebox(0,0)[lb]{\smash{\SetFigFont{12}{14.4}{\rmdefault}{\mddefault}{\updefault}
$a$
}}}
\put(7951,-6436){\makebox(0,0)[lb]{\smash{\SetFigFont{12}{14.4}{\rmdefault}{\mddefault}{\updefault}
$\overline{b_2}$
}}}
\put(7126,-5461){\makebox(0,0)[lb]{\smash{\SetFigFont{12}{14.4}{\rmdefault}{\mddefault}{\updefault}
$\overline{a}$
}}}
\put(11776,-4336){\makebox(0,0)[lb]{\smash{\SetFigFont{12}{14.4}{\rmdefault}{\mddefault}{\updefault}
$0$
}}}
\put(11326,-2986){\makebox(0,0)[lb]{\smash{\SetFigFont{12}{14.4}{\rmdefault}{\mddefault}{\updefault}
$c$
}}}
\put(10501,-2161){\makebox(0,0)[lb]{\smash{\SetFigFont{12}{14.4}{\rmdefault}{\mddefault}{\updefault}
$b_1$
}}}
\put(7951,-2236){\makebox(0,0)[lb]{\smash{\SetFigFont{12}{14.4}{\rmdefault}{\mddefault}{\updefault}
$b_2$
}}}
\end{picture}

\vspace{.1in}

Figure 1. $d\equiv 0\,(\mod 4)$, the sequence of zeros infinite in both
directions.
\end{center}
\newpage
\begin{center}
\begin{picture}(0,0)%
\includegraphics{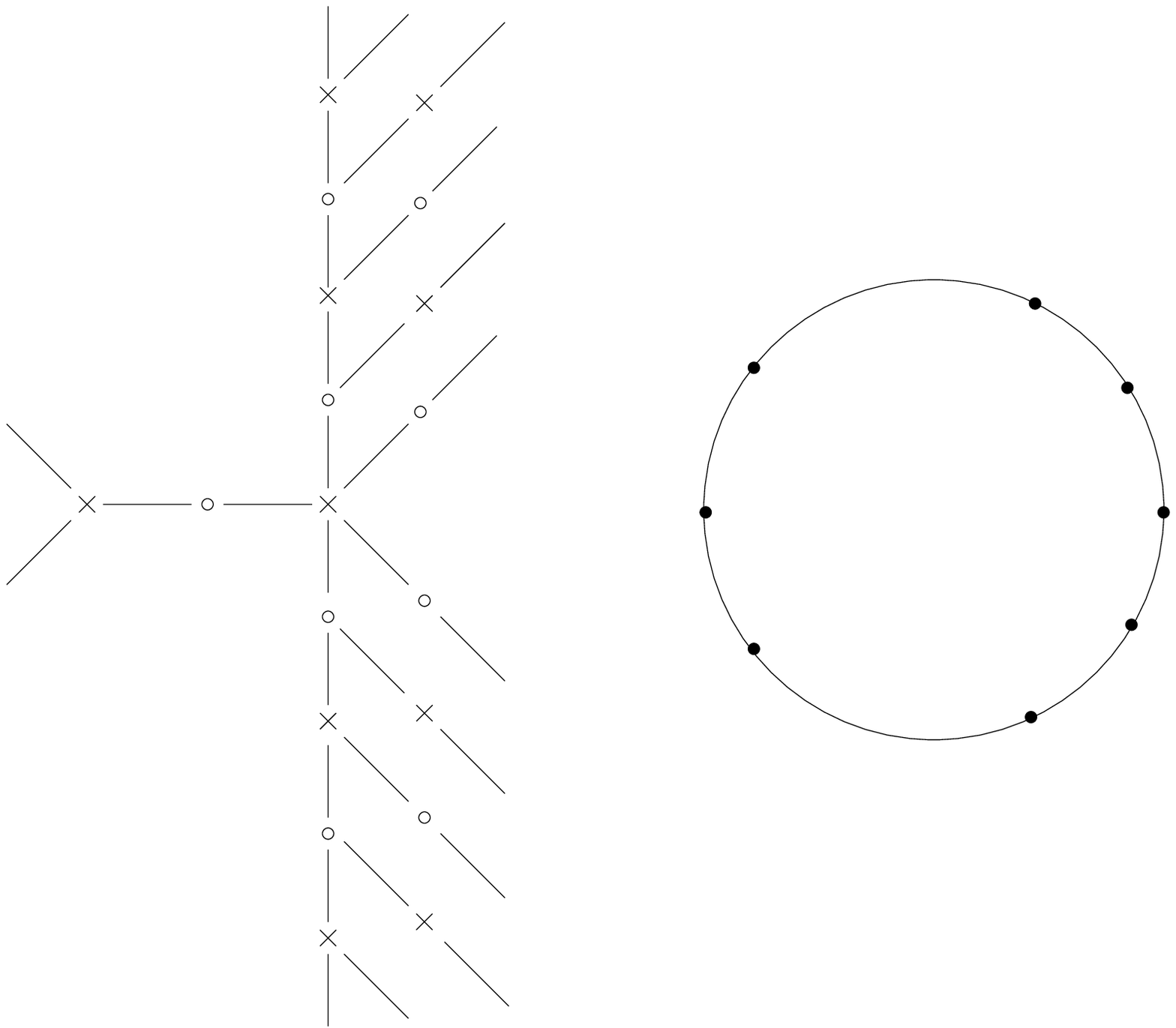}%
\end{picture}%
\setlength{\unitlength}{2763sp}%
\begingroup\makeatletter\ifx\SetFigFont\undefined%
\gdef\SetFigFont#1#2#3#4#5{%
  \reset@font\fontsize{#1}{#2pt}%
  \fontfamily{#3}\fontseries{#4}\fontshape{#5}%
  \selectfont}%
\fi\endgroup%
\begin{picture}(11025,9549)(676,-8848)
\put(10276,-1936){\makebox(0,0)[lb]{\smash{\SetFigFont{12}{14.4}{\rmdefault}{\mddefault}{\updefault}
$c_1$
}}}
\put(11326,-2761){\makebox(0,0)[lb]{\smash{\SetFigFont{12}{14.4}{\rmdefault}{\mddefault}{\updefault}
$c_2$
}}}
\put(11701,-4111){\makebox(0,0)[lb]{\smash{\SetFigFont{12}{14.4}{\rmdefault}{\mddefault}{\updefault}
$b$
}}}
\put(11251,-5386){\makebox(0,0)[lb]{\smash{\SetFigFont{12}{14.4}{\rmdefault}{\mddefault}{\updefault}
$\overline{c_2}$
}}}
\put(10276,-6286){\makebox(0,0)[lb]{\smash{\SetFigFont{12}{14.4}{\rmdefault}{\mddefault}{\updefault}
$\overline{c_1}$
}}}
\put(7501,-5686){\makebox(0,0)[lb]{\smash{\SetFigFont{12}{14.4}{\rmdefault}{\mddefault}{\updefault}
$\overline{a}$
}}}
\put(7201,-2686){\makebox(0,0)[lb]{\smash{\SetFigFont{12}{14.4}{\rmdefault}{\mddefault}{\updefault}
$a$
}}}
\put(3976,-2236){\makebox(0,0)[lb]{\smash{\SetFigFont{12}{14.4}{\rmdefault}{\mddefault}{\updefault}
$c_1$
}}}
\put(3976,-3286){\makebox(0,0)[lb]{\smash{\SetFigFont{12}{14.4}{\rmdefault}{\mddefault}{\updefault}
$0$
}}}
\put(3976,-1336){\makebox(0,0)[lb]{\smash{\SetFigFont{12}{14.4}{\rmdefault}{\mddefault}{\updefault}
$0$
}}}
\put(3976,-436){\makebox(0,0)[lb]{\smash{\SetFigFont{12}{14.4}{\rmdefault}{\mddefault}{\updefault}
$c_2$
}}}
\put(3901,389){\makebox(0,0)[lb]{\smash{\SetFigFont{12}{14.4}{\rmdefault}{\mddefault}{\updefault}
$0$
}}}
\put(4276,-4036){\makebox(0,0)[lb]{\smash{\SetFigFont{12}{14.4}{\rmdefault}{\mddefault}{\updefault}
$b$
}}}
\put(3976,-4936){\makebox(0,0)[lb]{\smash{\SetFigFont{12}{14.4}{\rmdefault}{\mddefault}{\updefault}
$0$
}}}
\put(3976,-5911){\makebox(0,0)[lb]{\smash{\SetFigFont{12}{14.4}{\rmdefault}{\mddefault}{\updefault}
$\overline{c_1}$
}}}
\put(3976,-7936){\makebox(0,0)[lb]{\smash{\SetFigFont{12}{14.4}{\rmdefault}{\mddefault}{\updefault}
$\overline{c_2}$
}}}
\put(3976,-6886){\makebox(0,0)[lb]{\smash{\SetFigFont{12}{14.4}{\rmdefault}{\mddefault}{\updefault}
$0$
}}}
\put(3901,-8836){\makebox(0,0)[lb]{\smash{\SetFigFont{12}{14.4}{\rmdefault}{\mddefault}{\updefault}
$0$
}}}
\put(1726,-5086){\makebox(0,0)[lb]{\smash{\SetFigFont{12}{14.4}{\rmdefault}{\mddefault}{\updefault}
$\overline{a}$
}}}
\put(1951,-3136){\makebox(0,0)[lb]{\smash{\SetFigFont{12}{14.4}{\rmdefault}{\mddefault}{\updefault}
$a$
}}}
\put(676,-4036){\makebox(0,0)[lb]{\smash{\SetFigFont{12}{14.4}{\rmdefault}{\mddefault}{\updefault}
$0$
}}}
\put(6751,-4111){\makebox(0,0)[lb]{\smash{\SetFigFont{12}{14.4}{\rmdefault}{\mddefault}{\updefault}
$0$
}}}
\end{picture}

\vspace{.1in}

Figure 2. $d\equiv 0\,(\mod 4)$, the sequence of zeros is finite.
\end{center}
\newpage
\begin{center}
\begin{picture}(0,0)%
\includegraphics{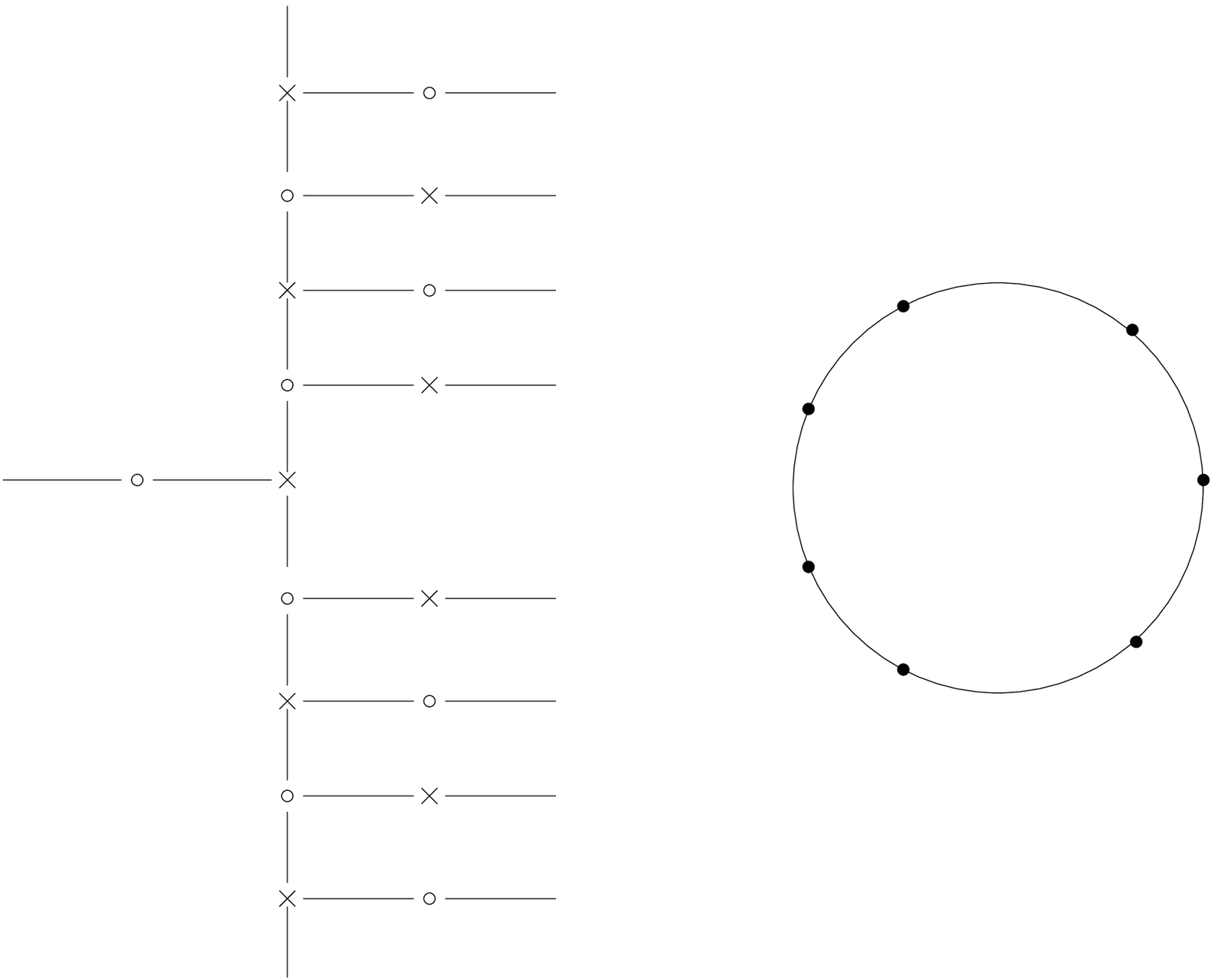}%
\end{picture}%
\setlength{\unitlength}{2763sp}%
\begingroup\makeatletter\ifx\SetFigFont\undefined%
\gdef\SetFigFont#1#2#3#4#5{%
  \reset@font\fontsize{#1}{#2pt}%
  \fontfamily{#3}\fontseries{#4}\fontshape{#5}%
  \selectfont}%
\fi\endgroup%
\begin{picture}(11562,9249)(589,-8698)
\put(4051, 14){\makebox(0,0)[lb]{\smash{\SetFigFont{12}{14.4}{\rmdefault}{\mddefault}{\itdefault}
$0$
}}}
\put(7651,-3136){\makebox(0,0)[lb]{\smash{\SetFigFont{12}{14.4}{\rmdefault}{\mddefault}{\itdefault}
$b_2$
}}}
\put(8701,-2161){\makebox(0,0)[lb]{\smash{\SetFigFont{12}{14.4}{\rmdefault}{\mddefault}{\itdefault}
$b_1$
}}}
\put(11251,-2386){\makebox(0,0)[lb]{\smash{\SetFigFont{12}{14.4}{\rmdefault}{\mddefault}{\itdefault}
$a$
}}}
\put(12151,-4036){\makebox(0,0)[lb]{\smash{\SetFigFont{12}{14.4}{\rmdefault}{\mddefault}{\itdefault}
$0$
}}}
\put(11251,-5836){\makebox(0,0)[lb]{\smash{\SetFigFont{12}{14.4}{\rmdefault}{\mddefault}{\itdefault}
$\overline{a}$
}}}
\put(8851,-6061){\makebox(0,0)[lb]{\smash{\SetFigFont{12}{14.4}{\rmdefault}{\mddefault}{\itdefault}
$\overline{b_1}$
}}}
\put(7951,-5086){\makebox(0,0)[lb]{\smash{\SetFigFont{12}{14.4}{\rmdefault}{\mddefault}{\itdefault}
$\overline{b_2}$
}}}
\put(1276,-5161){\makebox(0,0)[lb]{\smash{\SetFigFont{12}{14.4}{\rmdefault}{\mddefault}{\itdefault}
$\overline{a}$
}}}
\put(1276,-3061){\makebox(0,0)[lb]{\smash{\SetFigFont{12}{14.4}{\rmdefault}{\mddefault}{\itdefault}
$a$
}}}
\put(4051,-886){\makebox(0,0)[lb]{\smash{\SetFigFont{12}{14.4}{\rmdefault}{\mddefault}{\itdefault}
$b_2$
}}}
\put(4051,-1786){\makebox(0,0)[lb]{\smash{\SetFigFont{12}{14.4}{\rmdefault}{\mddefault}{\itdefault}
$0$
}}}
\put(4051,-2761){\makebox(0,0)[lb]{\smash{\SetFigFont{12}{14.4}{\rmdefault}{\mddefault}{\itdefault}
$b_1$
}}}
\put(4126,-4036){\makebox(0,0)[lb]{\smash{\SetFigFont{12}{14.4}{\rmdefault}{\mddefault}{\itdefault}
$0$
}}}
\put(4126,-5686){\makebox(0,0)[lb]{\smash{\SetFigFont{12}{14.4}{\rmdefault}{\mddefault}{\itdefault}
$\overline{b_1}$
}}}
\put(4126,-6586){\makebox(0,0)[lb]{\smash{\SetFigFont{12}{14.4}{\rmdefault}{\mddefault}{\itdefault}
$0$
}}}
\put(4126,-7561){\makebox(0,0)[lb]{\smash{\SetFigFont{12}{14.4}{\rmdefault}{\mddefault}{\itdefault}
$\overline{b_2}$
}}}
\put(4126,-8536){\makebox(0,0)[lb]{\smash{\SetFigFont{12}{14.4}{\rmdefault}{\mddefault}{\itdefault}
$0$
}}}
\end{picture}

\vspace{.1in}

Figure 3. $d\equiv 1\,(\mod 4)$.
\end{center}
\newpage
\begin{center}
\begin{picture}(0,0)%
\includegraphics{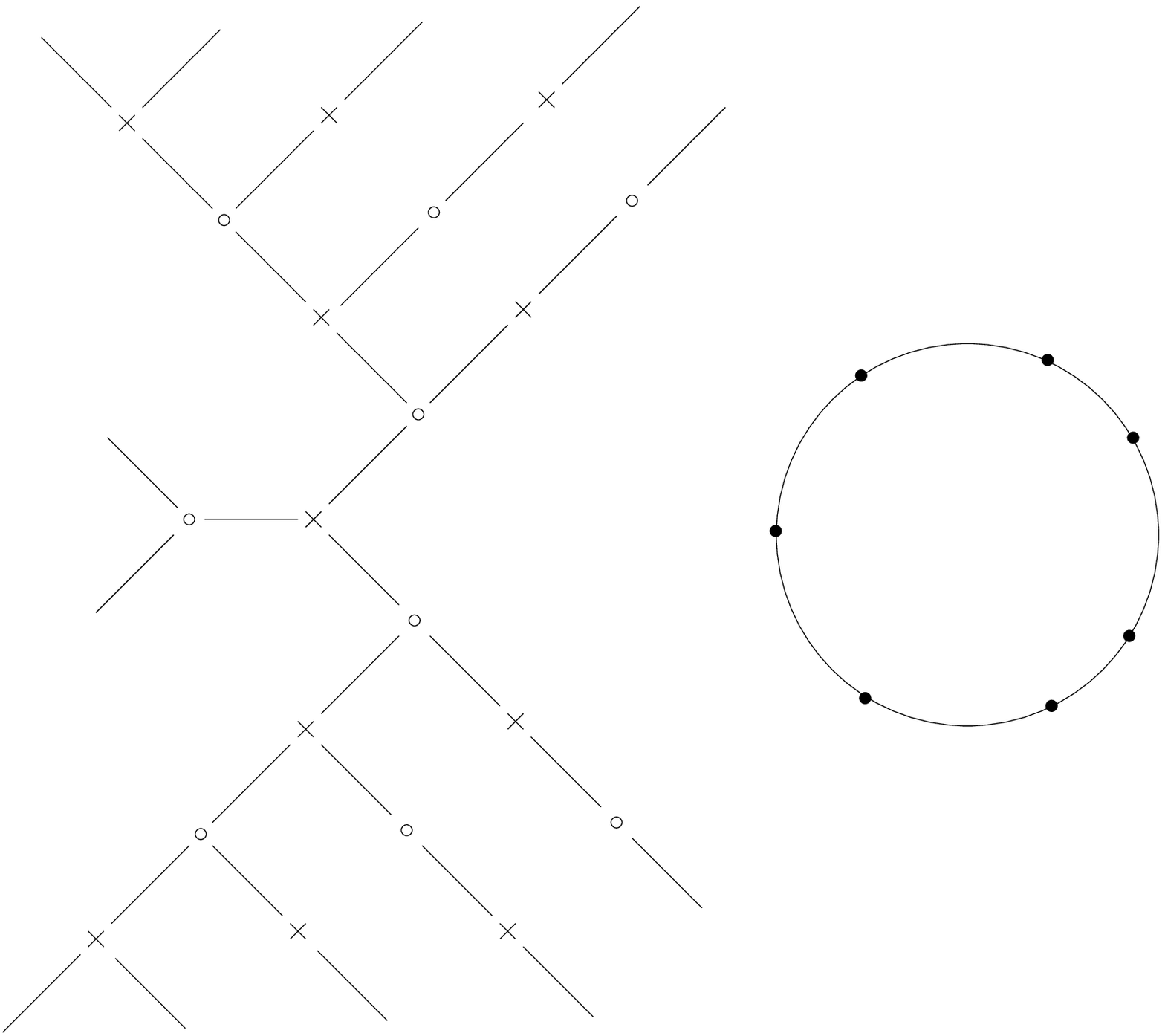}%
\end{picture}%
\setlength{\unitlength}{2763sp}%
\begingroup\makeatletter\ifx\SetFigFont\undefined%
\gdef\SetFigFont#1#2#3#4#5{%
  \reset@font\fontsize{#1}{#2pt}%
  \fontfamily{#3}\fontseries{#4}\fontshape{#5}%
  \selectfont}%
\fi\endgroup%
\begin{picture}(11178,9924)(-461,-8923)
\put(1576,-3061){\makebox(0,0)[lb]{\smash{\SetFigFont{12}{14.4}{\rmdefault}{\mddefault}{\updefault}
$a$
}}}
\put(526,314){\makebox(0,0)[lb]{\smash{\SetFigFont{12}{14.4}{\rmdefault}{\mddefault}{\updefault}
$0$
}}}
\put(6451,-4186){\makebox(0,0)[lb]{\smash{\SetFigFont{12}{14.4}{\rmdefault}{\mddefault}{\updefault}
$0$
}}}
\put(7351,-2536){\makebox(0,0)[lb]{\smash{\SetFigFont{12}{14.4}{\rmdefault}{\mddefault}{\updefault}
$a$
}}}
\put(9826,-2311){\makebox(0,0)[lb]{\smash{\SetFigFont{12}{14.4}{\rmdefault}{\mddefault}{\updefault}
$b_1$
}}}
\put(10576,-3136){\makebox(0,0)[lb]{\smash{\SetFigFont{12}{14.4}{\rmdefault}{\mddefault}{\updefault}
$b_2$
}}}
\put(10501,-5386){\makebox(0,0)[lb]{\smash{\SetFigFont{12}{14.4}{\rmdefault}{\mddefault}{\updefault}
$\overline{b_2}$
}}}
\put(9601,-6136){\makebox(0,0)[lb]{\smash{\SetFigFont{12}{14.4}{\rmdefault}{\mddefault}{\updefault}
$\overline{b_1}$
}}}
\put(1426,-361){\makebox(0,0)[lb]{\smash{\SetFigFont{12}{14.4}{\rmdefault}{\mddefault}{\updefault}
$b_2$
}}}
\put(3601,-4036){\makebox(0,0)[lb]{\smash{\SetFigFont{12}{14.4}{\rmdefault}{\mddefault}{\updefault}
$0$
}}}
\put(1351,-5011){\makebox(0,0)[lb]{\smash{\SetFigFont{12}{14.4}{\rmdefault}{\mddefault}{\updefault}
$\overline{a}$
}}}
\put(376,-4036){\makebox(0,0)[lb]{\smash{\SetFigFont{12}{14.4}{\rmdefault}{\mddefault}{\updefault}
$0$
}}}
\put(3151,-6211){\makebox(0,0)[lb]{\smash{\SetFigFont{12}{14.4}{\rmdefault}{\mddefault}{\updefault}
$\overline{b_1}$
}}}
\put(2326,-7186){\makebox(0,0)[lb]{\smash{\SetFigFont{12}{14.4}{\rmdefault}{\mddefault}{\updefault}
$0$
}}}
\put(226,-8911){\makebox(0,0)[lb]{\smash{\SetFigFont{12}{14.4}{\rmdefault}{\mddefault}{\updefault}
$0$
}}}
\put(1276,-8086){\makebox(0,0)[lb]{\smash{\SetFigFont{12}{14.4}{\rmdefault}{\mddefault}{\updefault}
$\overline{b_2}$
}}}
\put(3301,-2161){\makebox(0,0)[lb]{\smash{\SetFigFont{12}{14.4}{\rmdefault}{\mddefault}{\updefault}
$b_1$
}}}
\put(2476,-1261){\makebox(0,0)[lb]{\smash{\SetFigFont{12}{14.4}{\rmdefault}{\mddefault}{\updefault}
$0$
}}}
\put(7651,-6061){\makebox(0,0)[lb]{\smash{\SetFigFont{12}{14.4}{\rmdefault}{\mddefault}{\updefault}
$\overline{a}$
}}}
\end{picture}

\vspace{.1in}

Figure 4. $d\equiv 2\,(\mod 4).$
\end{center}
\newpage
\begin{center}
\begin{picture}(0,0)%
\includegraphics{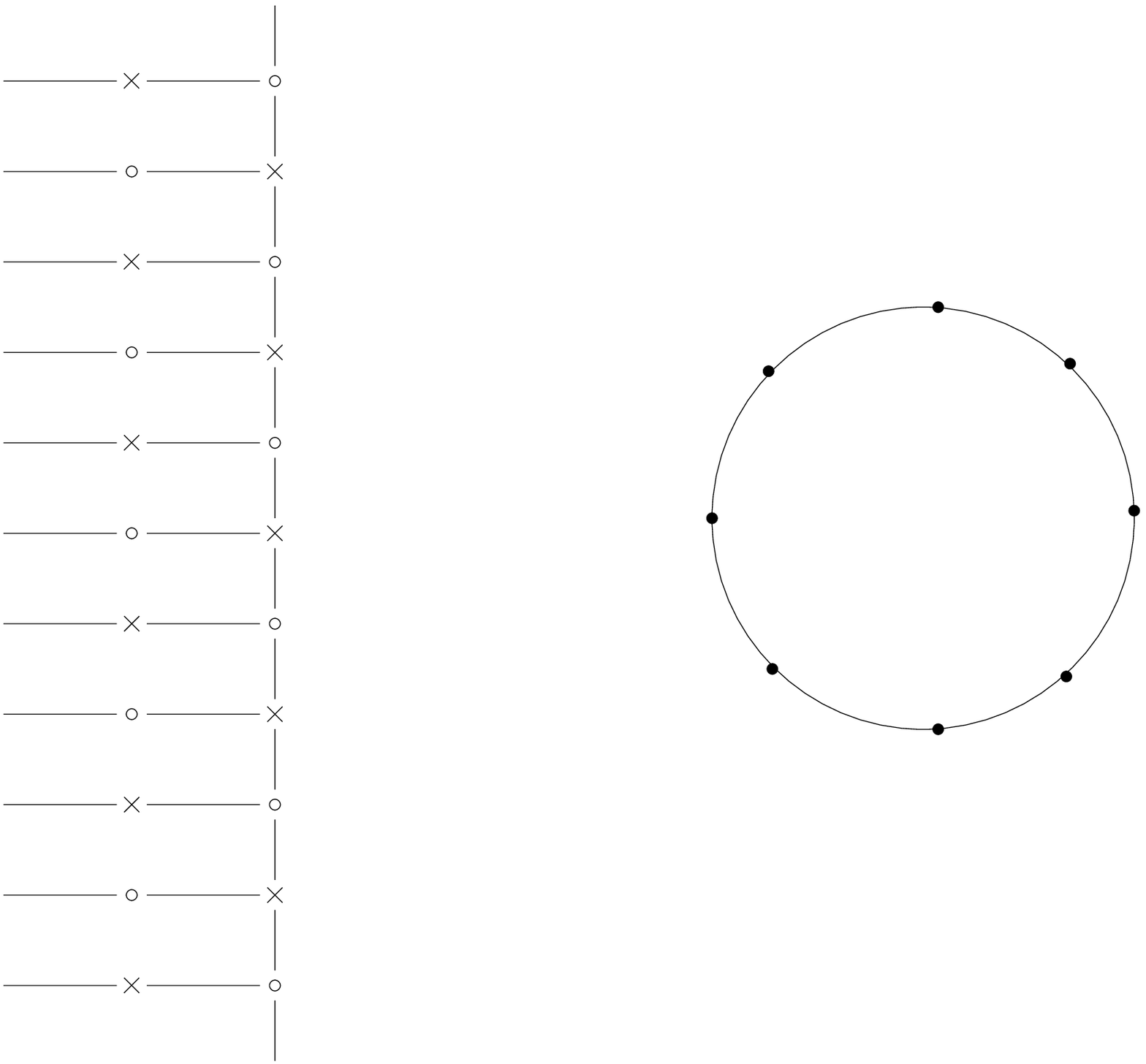}%
\end{picture}%
\setlength{\unitlength}{2763sp}%
\begingroup\makeatletter\ifx\SetFigFont\undefined%
\gdef\SetFigFont#1#2#3#4#5{%
  \reset@font\fontsize{#1}{#2pt}%
  \fontfamily{#3}\fontseries{#4}\fontshape{#5}%
  \selectfont}%
\fi\endgroup%
\begin{picture}(11412,10524)(739,-9523)
\put(9751,-1861){\makebox(0,0)[lb]{\smash{\SetFigFont{12}{14.4}{\rmdefault}{\mddefault}{\updefault}
$b_1$
}}}
\put(11251,-6061){\makebox(0,0)[lb]{\smash{\SetFigFont{12}{14.4}{\rmdefault}{\mddefault}{\updefault}
$\overline{b_2}$
}}}
\put(9826,-6511){\makebox(0,0)[lb]{\smash{\SetFigFont{12}{14.4}{\rmdefault}{\mddefault}{\updefault}
$\overline{b_1}$
}}}
\put(8101,-5911){\makebox(0,0)[lb]{\smash{\SetFigFont{12}{14.4}{\rmdefault}{\mddefault}{\updefault}
$\overline{a}$
}}}
\put(7276,-4186){\makebox(0,0)[lb]{\smash{\SetFigFont{12}{14.4}{\rmdefault}{\mddefault}{\updefault}
$0$
}}}
\put(7876,-2611){\makebox(0,0)[lb]{\smash{\SetFigFont{12}{14.4}{\rmdefault}{\mddefault}{\updefault}
$a$
}}}
\put(11401,-2386){\makebox(0,0)[lb]{\smash{\SetFigFont{12}{14.4}{\rmdefault}{\mddefault}{\updefault}
$b_2$
}}}
\put(12151,-4111){\makebox(0,0)[lb]{\smash{\SetFigFont{12}{14.4}{\rmdefault}{\mddefault}{\updefault}
$b$
}}}
\put(1426,-6586){\makebox(0,0)[lb]{\smash{\SetFigFont{12}{14.4}{\rmdefault}{\mddefault}{\updefault}
$\overline{b_1}$
}}}
\put(1501,-5686){\makebox(0,0)[lb]{\smash{\SetFigFont{12}{14.4}{\rmdefault}{\mddefault}{\updefault}
$0$
}}}
\put(1426,-7486){\makebox(0,0)[lb]{\smash{\SetFigFont{12}{14.4}{\rmdefault}{\mddefault}{\updefault}
$0$
}}}
\put(1426,-8386){\makebox(0,0)[lb]{\smash{\SetFigFont{12}{14.4}{\rmdefault}{\mddefault}{\updefault}
$\overline{b_2}$
}}}
\put(1426,-9136){\makebox(0,0)[lb]{\smash{\SetFigFont{12}{14.4}{\rmdefault}{\mddefault}{\updefault}
$0$
}}}
\put(4426,-4261){\makebox(0,0)[lb]{\smash{\SetFigFont{12}{14.4}{\rmdefault}{\mddefault}{\updefault}
$b$
}}}
\put(1501,-4786){\makebox(0,0)[lb]{\smash{\SetFigFont{12}{14.4}{\rmdefault}{\mddefault}{\updefault}
$\overline{a}$
}}}
\put(1501,-3811){\makebox(0,0)[lb]{\smash{\SetFigFont{12}{14.4}{\rmdefault}{\mddefault}{\updefault}
$a$
}}}
\put(1501,-2911){\makebox(0,0)[lb]{\smash{\SetFigFont{12}{14.4}{\rmdefault}{\mddefault}{\updefault}
$0$
}}}
\put(1501,-2161){\makebox(0,0)[lb]{\smash{\SetFigFont{12}{14.4}{\rmdefault}{\mddefault}{\updefault}
$b_1$
}}}
\put(1501,-1186){\makebox(0,0)[lb]{\smash{\SetFigFont{12}{14.4}{\rmdefault}{\mddefault}{\updefault}
$0$
}}}
\put(1501,-286){\makebox(0,0)[lb]{\smash{\SetFigFont{12}{14.4}{\rmdefault}{\mddefault}{\updefault}
$b_2$
}}}
\put(1501,464){\makebox(0,0)[lb]{\smash{\SetFigFont{12}{14.4}{\rmdefault}{\mddefault}{\updefault}
$0$
}}}
\end{picture}

\vspace{.1in}

Figure 5. $d\equiv 3\,(\mod 4).$
\end{center}

{\em 

Purdue University

West Lafayette IN

eremenko@math.purdue.edu
\vspace{.2in}

University of Michigan

Ann Arbor MI

merenkov@umich.edu}

\end{document}